\documentclass[final, nomarks]{dmtcs-episciences}
\usepackage[latin1]{inputenc}
\usepackage{amsmath,amsthm,amsfonts,amssymb}
\usepackage{color}
\usepackage{epsfig}
\usepackage{mathrsfs}
\usepackage{subfigure}
\usepackage{graphicx}
\usepackage[ruled,vlined,linesnumbered]{algorithm2e}
\def\slika #1{\begin{center} \epsffile{#1} \end{center}}
\newtheorem{theorem}{Theorem}

\newtheorem{lemma}{Lemma}

\newtheorem{corollary}{Corollary}
\newtheorem{remark}{Remark}

\author[Arun Anil, Manoj Changat]{Arun Anil\thanks{arunanil93@gmail.com}
	\and Manoj Changat\thanks{mchangat@keralauniversity.ac.in}}

\title[Recognition of chordal C-I graphs and C-I cographs]{Recognition of chordal graphs and cographs which are Cover-Incomparability graphs}

\affiliation{Department of Futures Studies, University of Kerala, Thiruvananthapuram, Kerala, India.}

\begin{document}
\publicationdata{vol. 26:3}{2024}{12}{10.46298/dmtcs.11657}{2023-07-27; 2023-07-27; 2024-04-11; 2024-09-23}{2024-10-02}
\maketitle	
	\begin{abstract}
		Cover-Incomparability graphs (C-I graphs) are an interesting class of graphs from posets. A C-I graph is a graph from a poset $P=(V,\le)$ with vertex set $V$, and the edge-set is the union of edge sets of the cover graph and the incomparability graph of the poset. The recognition of the C-I graphs is known to be NP-complete (Maxov\'{a} et al., Order 26(3), 229--236(2009)).
		In this paper, we prove that chordal graphs having at most two independent simplicial vertices are exactly the chordal graphs which are also C-I graphs. A similar result is obtained for cographs as well. 
		Using the structural results of these graphs, we derive linear-time recognition algorithms for chordal graphs and cographs which are C-I graphs.
	\end{abstract}
	
	\section{Introduction}\label{section_introduction}
	
  The cover-incomparability graphs of posets, or shortly C-I graphs form an interesting class of graphs from posets. These graphs were introduced in \cite{bckkmm-07} as the underlying graphs of the so-called standard transit function of posets. The C-I graphs are precisely the graphs whose edge set is the union of edge sets of the cover graph and the incomparability graph (complement of a comparability graph) of a poset. Graph-theoretic characterization of C-I graphs has not been known until now and, moreover, the recognition complexity of C-I graphs is NP-complete (Maxov\'{a} et al. \cite{mpt-09}). Hence, problems in C-I graphs are focused on identifying the structure and characterization of well-known graph families, which are C-I graphs. Such C-I graphs studied include the family of split graphs, block graphs \cite{bcgmm-09}, cographs \cite{cograph}, Ptolemaic graphs \cite{matu-13}, distance-hereditary graphs \cite{matu-13}, and $k$-trees \cite{mdp-14}. The C-I graphs were identified among the planar and chordal graphs along with new characterizations of the Ptolemaic graphs, respectively, in \cite{amtb-01} and \cite{am-02}. It is also interesting to note that every C-I graph has a Ptolemaic C-I graph as a spanning subgraph \cite{amtb-01}. C-I graphs are also studied among the comparability graphs \cite{am-03}. The effect of composition operation, lexicographic, and strong products of C-I graphs was studied in a recent paper \cite{am-04}.

 Another approach in the theory of C-I graphs is to study posets whose cover-incomparability graphs exhibit specific properties \cite{bcgmm-09, bckkmm-07, cograph}. Characterizations of posets whose cover-incomparability graphs are distance-hereditary and Ptolemaic were given in \cite{bckkmm-07}, while those for cographs appear in \cite{cograph}, and for block graphs and split graphs in \cite{bcgmm-09}, in terms of forbidden isometric subposets. However, the characterization in terms of forbidden `isometric subposets' in \cite{bcgmm-09, bckkmm-07, cograph} was an error, and it was resolved in \cite{janbok-02} by introducing $\lhd$-preserving subposets instead of isometric subposets.
 
Chordal graphs or triangulated graphs are graphs that do not have induced cycles of length greater than three. Together with cographs, which are exactly $P_4$-free graphs, they  form two well-studied graph classes that have wider applications beyond mathematics. These two classes of graphs and their subclasses are also studied in many contexts, including theoretical, practical, and algorithmic interests.  A vertex $v$ of a graph $G$ is a simplicial vertex if $v$ together with all its adjacent vertices forms a clique (a complete subgraph) in $G$. %Simplicial elimination play a crucial role in chordal graphs. 
 It is a well-known fact that every chordal graph contains at least one simplicial vertex, and a non complete chordal graph contains at least two nonadjacent simplicial vertices. Chordal graphs that contain exactly two non-adjacent simplicial vertices form a special class of chordal graphs. %By a chordal cograph, we mean a graph which are both a chordal graph and a cograph.
	
 %In a similar manner, chordal cographs that contain exactly two non-adjacent simplicial vertices are a proper subclass of the special class of chordal graphs mentioned above. Also, the join of such chordal cographs form a special class of cographs.\\
 
The recognition algorithm for an arbitrary chordal graph \cite{trajan,rose_chordal-2} and a cograph \cite{corneil} is well known and has linear-time complexity.  In this paper, we identify chordal graphs having exactly two non-adjacent simplicial vertices and prove that they are C-I graphs. In other words, we prove that the intersection of the class of chordal graphs and the class of C-I graphs is precisely the class of chordal graphs having exactly two non-adjacent simplicial vertices. We further develop a linear-time algorithm for recognizing such chordal graphs. It is interesting to observe, from \cite{amtb-01, bcgmm-09}, that in the class of C-I graphs, the classes of chordal, strongly chordal, and interval graphs are one and the same.

In this paper, similar to the chordal C-I graphs, we develop a linear-time recognition algorithm for cographs which are also C-I graphs from the structural results already proved in \cite{cograph}. For this algorithm, we make use of the special type of cotree structure of cographs.

	Chordal graphs are precisely the intersection graphs of subtrees of a tree \cite{gravril}. Most of the information contained in a chordal graph is captured in its clique tree representation, which is useful for its algorithmic applications \cite{gravril, shibata}. 
	
	Cographs are exactly the $P_4$-free graphs and the class of cographs has been intensively studied since its definition by Seinsche \cite{seinsche}.  The cographs appear as comparability graphs of series-parallel partial orders \cite{jung}, and can be generated from the single-vertex graph $K_1$ by complementation and disjoint union operations. It is well known that any cograph has a canonical tree representation called a cotree. This tree decomposition scheme for cographs is a particular case of modular decomposition \cite{Gallai} that applies to arbitrary graphs. Indeed, the algorithm that computes the modular tree decomposition of an arbitrary graph in linear-time can also recognize cographs in linear-time.  In 1994, linear-time modular decomposition algorithms were designed independently by Cournier and Habib \cite{habib-c-1} and by McConnell and Spinrad \cite{McConnell-1994}. In 2001, Dahlhaus et al. \cite{ Dahlhaus-2001} proposed a simpler algorithm. %Unfortunately, because they build the decomposition tree, \textcolor{red}{all these algorithms involve to maintain complicated data structures}. 
    In 2004, Habib and Paul \cite{habib-2005} proposed a new algorithm, which is not incremental, and instead of building the cotree directly, it first computes a special ordering of the vertices, namely a factorizing permutation, using a  very efficient partition refinement techniques via two elementary refinement rules.

	A search of a graph visits all vertices and edges of the graph and will visit a new vertex only if it is adjacent to some previously visited vertex. The two fundamental search strategies are Breadth-First Search (BFS) and Depth-First Search (DFS). As the names indicate, BFS visits all previously unvisited neighbors of the currently visited vertex before visiting the previously unvisited non-neighbors. 
	Several greedy recognition algorithms for chordal graphs are known. The most famous is Lex-BFS \cite{rose_chordal-2}, a variant of $BFS$, introduced by Rose et al. in \cite{rose_chordal-2} and Maximum Cardinality Search (MCS for short) \cite{trajan}. Both algorithms are linear.  The recognition of chordal graphs involves two distinct phases: the execution of MCS or Lex-BFS in order to compute an elimination ordering and a checking procedure to decide whether this elimination ordering is perfect (PEO). By employing Lex-BFS as a basic method to check the chordality and to find the perfect elimination ordering, the recognition algorithm for the chordal graphs can be done in linear-time for chordal graphs which are  C-I graphs. For cographs, using the BFS, we can check whether a given rooted tree is the cotree of the C-I cograph in linear-time because of the special structure of the cotrees of these graphs.

	The rest of this section is organized as follows. In Section ~\ref{section_Prelim}, we begin with some preliminaries. In Section~\ref{section_chordal}, we characterize chordal C-I graphs in terms of the number of independent simplicial vertices and present a linear-time algorithm for recognizing chordal C-I graphs based on this characterization. In Section~\ref{Section_cograph}, we present a structure of the C-I cograph and characterize the structure of the cotree of C-I cographs and, using this, present a linear-time recognition algorithm for C-I cographs.
	
	\section{Preliminaries}\label{section_Prelim}
	%Some preliminary definitions and results that were used in this paper are discussed in this section.
	
	A \textit{partially ordered set} or \emph{poset} $P=(V,\leq)$ consists of a nonempty set $V$ and a reflexive, antisymmetric, transitive relation $\leq$ on $V$, denoted as $P= (V, \le)$. We call $u \in V$ an element of $P$. If $u\leq v$ or $v\leq u$ in $P$, we say $u$ and $v$ are \textit{comparable}, otherwise \textit{incomparable}. If $u\leq v$ but $u\neq v$, then we write $u<v$.  If $u$ and $v$ are in $V$, then $v$ \textit{covers} $u$ in $P$ if $u<v$ and there is no $w$ in $V$ with $u<w<v$, denoted by $u \lhd v$. We write $u \lhd\lhd v$ if $u < v$ but not $u \lhd v$. By $u||v$, we mean that $u$ and $v$ are incomparable elements of $P$. In this paper, we consider only posets defined on finite sets. 
	Let $V' \subseteq V$ and $Q = (V',\leq ')$ be a poset, $Q$ is called a {\em subposet} of $P$, if $u \leq ' v$ if and only if $u \leq v$, for any $u,v \in V'$. The subposet $Q=(V',\leq)$ is a {\em chain (antichain)} in $P$, if every pair of elements of $V'$ is comparable (incomparable) in $P$. A chain of maximum cardinality is named as the {\em height} of $P$ denoted as $h(P)$. An element $u$ in $P$ is a \emph{minimal (maximal)} if there is no $x\in V$ such that $x\le u (u\le x)$ in $P$. 	A finite {\em ranked poset} (also known as
	\emph{graded poset}~\cite{bw-00}) is a
	poset $P =(V, \le )$ that is equipped with a rank function $\rho:
	V\rightarrow \mathbb{Z}$ satisfying:
	\begin{itemize}
		\item $\rho$ has value $0$ on all minimal elements of $P$, and 
		\item $\rho$ preserves covering relations: if $a\lhd b$ then $\rho (b)=\rho (a) + 1$.
	\end{itemize}
	A ranked poset $P$ is said to be {\em complete} if for every $i$, every element of rank $i$ covers all elements of rank $i-1$. For a completely ranked poset $P=(V, \leq)$ we say that the element $v \in V$ is at height $i$ if $\rho(v)=i-1$. For disjoint posets, $P_1 = (V_1,\le)$ and
$P_2 =(V_2,\le)$, the \emph{sum} of the posets, $Z = P_1+P_2$ is defined as the poset $(V_1\cup V_2,\le)$,
where $z_1\in z_2$ in $Z$, if and only if $z_1\le z_2$ in $P_1$ or $z_1\le z_2$ in $P_2$. We refer to \cite{bw-00}, for notions of posets.
	
	Let $G=(V, E)$ be a connected graph, vertex set and edge set of $G$ denoted as $V(G)$ and $E(G)$ respectively, the complement of $G$ is denoted as $\overline{G}$.  For a vertex $v\in V(G)$, the set of all vertices adjacent to $v$ is called the {\em open neighborhood} of $v$ and is denoted by $N(v)$. The set consisting of the open neighborhood and the vertex $v$ is the {\em closed neighborhood} of $v$ and is denoted by $N[v]$. A \emph{tree} is a graph in which two vertices are connected by exactly one path. Let $T$ be a rooted tree and two vertices $x$ and $y$ in $T$, we say that $x$ is an ancestor of $y$ and $y$ is a descendant of $x$ if $x$ lies on the path from $y$ to the root of $T$. For a set of leaves $S$ of $T$, we say that the \emph{lowest common ancestor (LCA)} of $S$ is the internal node $v$ of $T$ such that $v$ is the root of the smallest rooted subtree of $T$ containing $S$.
	A graph $H$ is said to be a \emph{subgraph} of $G$ if $V(H)\subseteq V(G)$ and $E(H)\subseteq E(G)$.  $H$ is an \emph{induced subgraph} of $G$ if for $u,v\in V(H)$ and $uv\in E(G)$ implies $uv\in E(H)$. A graph $G$ is said to be \textit{H-free} if $G$ has no induced subgraph isomorphic to $H$. A {\em complete graph}  is a graph whose vertices are pairwise adjacent, denoted $K_n$, a set $S \subseteq V(G)$ is a \emph{clique} if the subgraph of $G$ induced by $S$ is a complete graph, and a \emph{ maximum clique} is a clique that is not contained by any other clique. A vertex $v$ is called \emph{simplicial vertex} if its neighborhood induces a complete subgraph. An \emph{independent set} in a graph is a set of pairwise non-adjacent vertices. 
	If graphs $G_{1}$ and $G_{2}$ have disjoint vertex set $V_{1}$ and $V_{2}$ and edge set $E_{1}$ and $E_{2}$ respectively, then their \emph{union} $G=G_{1} \cup G_{2}$ has $V=V_{1} \cup V_{2}$ and $E=E_{1} \cup E_{2}$ and their \emph{join}, denoted by $G_{1}\vee G_{2}$, consists of $G=G_{1} \cup G_{2}$ and all edges joining $V_{1}$ with $V_{2}$.
	
	A graph $G$ is {\em chordal} if it contains no induced cycles of length more than 3. A graph $G$ is {\em distance-hereditary} if every induced path is also the shortest path in $G$.  A graph $G$ is {\em Ptolemaic} if it is distance-hereditary and chordal. Equivalently, $G$ is Ptolemaic if and only if it is a 3-fan-free chordal graph. $P_4$-free graphs are called \emph{cographs}. A graph $G$ that is both chordal and cograph is called \emph{ chordal cograph}, also known as \emph{trivially perfect graph}. 
 Let $P=(V,\le)$ be a poset. A graph $G$ is the \emph{cover graph} of  $P$ if  $V(G)=V$ and $uv\in E(G)$  if and only if either  $u\lhd v$ or $v\lhd u$ in $P$. Similarly, $G$ is the \emph{comparability graph} of $P$ if $V(G)=V$ and $uv\in E(G)$   if and only if $u$ and $v$ are comparable in $P$. The \emph{incomparability graph} is the complement of the comparability graph. Finally, the \emph{cover-incomparability graph (C-I graph)} of a poset $P= (V,\leq)$ denoted as $G_P$ is the graph $G=(V,E)$, where $uv\in E(G)$, if $u\lhd v$ or $v\lhd u$ or $u||v$ in $P$. A graph is a C-I graph if it is the C-I graph of some poset $P$. It may be noted that the class of C-I graphs is not a \textit{hereditary class} in the sense that every induced subgraph of a C-I graph need not be a C-I graph. We call a graph that is both a chordal and a C-I graph as \emph{ a chordal C-I graph}. Similarly, we use the term \emph{ Ptolemaic C-I graph, C-I cograph, chordal C-I cograph}, etc. to denote the C-I graph which is Ptolemaic, cograph, chordal cograph, etc.

  \begin{figure}[htb]	
		\epsfxsize=11truecm \slika{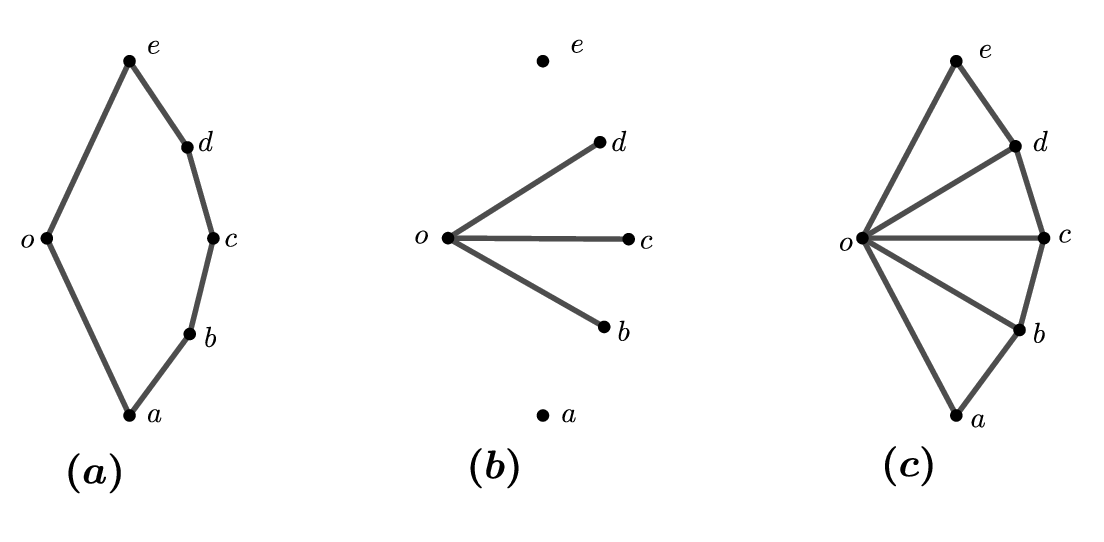}\vspace{-.85cm} \caption{(a) A poset $P$, (b) incomparability graph of the poset $P$, and (c) cover-incomparability graph of the poset $P$.}\label{fig:eg_C-I}
	\end{figure}
	In Figure~\ref{fig:eg_C-I} (a), we illustrates the Hasse diagram of a poset $P$, which is isomorphic to the cover graph of $P$. Figure~\ref{fig:eg_C-I}(b), depicts the incomparability graph of the poset $P$, and in Figure~\ref{fig:eg_C-I}(c), we  depicts the C-I graph of $P$. As already mentioned, the class of C-I graphs is not hereditary; for example, the 4-fan is a C-I graph (as shown in Figure~\ref{fig:eg_C-I}(c)), but the claw graph is an induced subgraph of this graph, which is not a C-I graph.
	
	Now we recall some basic properties of posets and their C-I graphs.
	
	\begin{lemma}\cite{bckkmm-07}\label{basic_result}
		Let P be a poset. Then
		\begin{enumerate}
			\item[(i)] the C-I graph of P is connected;
			\item[(ii)] points of P that are independent in the C-I graph of P lie on a common chain;
			\item[(iii)] an antichain of P corresponds to a complete subgraph in the C-I graph of P; 
			\item[(iv)] the C-I graph of P contains no induced cycles of length greater than 4.
		\end{enumerate}
	\end{lemma}
	\begin{lemma} \cite{matu-13} \label{CI3notsimplicial}
		If $G$ is a C-I graph, then $G$ does not contain 3 independent simplicial vertices.
	\end{lemma}

	\begin{lemma} \cite{matu-13} \label{max:mini} Let $P$ be a poset and $G$ its C-I graph. If $v$ is a simplicial vertex in $G$, then $v$ is a maximal or a minimal element of $P$.
	\end{lemma}
	
	\begin{theorem}\cite{mdp-14}\label{maxi_mini}
		Let $G=(V, E)$ be a C-I graph of a poset $P$ and $v\in V$ a minimal or maximal element in $P$. Then $G\backslash v$ is also a C-I graph.
	\end{theorem}

	In the following sections, we discuss algorithms for recognizing chordal C-I graphs and C-I cographs.
	\section{Structure of chordal C-I graphs}\label{section_chordal}
	
	In this section, we study the structure of chordal C-I graphs. For that, we prove that a chordal graph having exactly two independent vertices is a chordal graph that is also a C-I graph.
	
	If $v$ is a simplicial vertex of a chordal graph $G$, then the closed neighborhood $N[v]$ is a clique in $G$. So removal of $v$ does not affect the chordal property of $G\setminus \{v\}$. This is stated in the following remark.
	\begin{remark}\label{simpicial_minus_chordal}
		If $G$ is a chordal graph and $v$ is a simplicial vertex, then $G\setminus v$ is chordal.
	\end{remark}
	\begin{lemma}\label{G-v_simplicial}
		Let $G$ be a chordal graph, and $v$ be a simplicial vertex of $G$. If $G$ has exactly two independent simplicial vertices, then either $G\setminus v$ is chordal and has exactly two independent simplicial vertices or $G\setminus v$ is a complete graph.
	\end{lemma}
	\begin{proof}
		Let $G$ be a chordal graph with exactly two independent simplicial vertices and $v$ be a simplicial vertex. Then by Remark~\ref{simpicial_minus_chordal}, $G\setminus v$ is a chordal graph. Since $v$ is a simplicial vertex, the removal of $v$ from $G$ does not increase the number of independent simplicial vertices in $G\setminus v$. And if $G\setminus v$ is not a complete graph, then removal of $v$ from $G$ does not decrease the number of independent simplicial vertex in $G\setminus v$. That is, $G\setminus v$ has exactly two independent simplicial vertices, otherwise $G\setminus v$ is a complete graph. Hence the lemma.
	\end{proof}
	
	From Lemmas \ref{max:mini}-\ref{G-v_simplicial} and Theorem \ref{maxi_mini}, we get the following corollary.
	
	\begin{corollary}\label{coro_G-v_2_indpendent}
		Let $G$ be a chordal C-I graph and $v$ be a simplicial vertex of $G$. If $G$ has exactly two independent simplicial vertices, then $G\setminus v$ is a chordal C-I graph and has exactly two independent simplicial vertices or $G\setminus v$ is a complete graph.
	\end{corollary}
	
	A \emph{pendant clique} of a graph $G$ is a clique that contains a clique separator $C$ such that one of the components obtained after removing $C$ is a single vertex. Observe that if $G$ is a chordal graph, then $G$ has a pendant clique.\\

	Let $G$ be a chordal graph having at most two independent simplicial vertices. We denote by $\mathscr{G}$  the family of chordal graphs obtained from $G$  by adding a vertex $v$ to $G$ such that $v$ is adjacent to some or all vertices in a pendant clique of $G$. 
	
	The following remark is immediate. 
	
 \begin{remark} \label{lemma_adding_simplicial}
		Let $G'$ be a chordal graph having exactly two independent simplicial vertices. If $G$ is the graph obtained by adding a simplicial vertex $v$  to $G'$ such that the resulting graph has exactly two independent simplicial vertices, then $G$ belongs to the family $\mathscr{G}$.
	\end{remark}

	\begin{lemma}\label{G_k_chordal}
		Let $G'$ be a chordal C-I graph having exactly two independent simplicial vertices, and let $G$ be a graph obtained by adding a simplicial vertex $v$ to $G'$. If $G$ has exactly two independent simplicial vertices, then $G$ is a chordal C-I graph.
	\end{lemma}
	\begin{proof}
		Let $G$ be a chordal graph obtained by adding a simplicial vertex $v$ to a chordal $G'$ that has exactly two independent simplicial vertices such that $G$ also has exactly two independent simplicial vertices. Then by Remark~\ref{lemma_adding_simplicial}, the vertex $v$ is such that $v$ is adjacent to some or all vertices in a pendant clique in $G'$. 
		
		Now we need to prove that $G$ is a chordal C-I graph. Since $G'$ is both a C-I graph and chordal, let $P'$ be a poset such that $G_{P'} \cong G'$. Since the vertex $v$ is added to some or all vertices in a pendant clique in $G'$, the graph $G$ is chordal. Since $G_P'$ is chordal and has exactly two independent simplicial vertices, $G_P'$ has exactly two pendant cliques, and the two pendent cliques are formed, respectively, by $S$ and $S'$, where $S$ and $S'$ are defined as follows.\\
		$M=\{u\in P'\mid u \text{ is a maximal element of } P'\}$, 
		$S_1=\{w\in P'\mid w\lhd u, u\in M \}$, $M'=\{u\in P'\mid u \text{ is a minimal element of } P'\}$,
		$S_1'=\{w\in P'\mid u\lhd w, u \in M'\}$ 
		$S=M \cup S_1\setminus \{u\in S_1 \mid \exists w\in S_1, u\lhd\lhd w \text{ in }P'\}$ and \\ $S'= M'\cup S_1'\setminus \{u\in S_1'\mid \exists w\in S_1', w\lhd\lhd u \text{ in }P'\}$. Now $v$ is adjacent to some or all vertices of the set $S$ or $S'$ in $G$. If $v$ is adjacent to some vertices of $S$, then $v$ must be adjacent to all the elements of $M$, as otherwise there will be more than two independent simplicial vertices. Similarly, if $v$ is adjacent to some vertices of $S'$, then $v$ must be adjacent to all elements of $M'$. Now, we construct a poset $P$ from $P'$ as follows.  
		
		If $v$ adjacent to some vertices of $S$:
		\begin{itemize}
			
			\item If $v$ is adjacent to only the elements in $M$, then $P$ is constructed from $P'$ with the covering relation defined as $u\lhd v$ for all $u\in M$.
			\item If $v$ is adjacent to the elements of $S_0$, where $M\subset S_0\subset S$, then $P$ is constructed from $P'$ with the covering relation as $u'\lhd v$ for all $u'\in S_0\setminus M$ and $u''\lhd u\lhd v$ for all $u''\in S\setminus S_0$ for some $u\in M$.
			\item If $v$ is adjacent to all elements of $S$, then $P$ is constructed from $P'$ with the covering relation defined as $u\lhd v$ for all $u \in S\setminus M$.
		\end{itemize}
		Similarly,  if $v$ adjacent to some vertices of $S'$:
		\begin{itemize}
			
			\item if $v$ is adjacent to only the elements in $M'$ then $P$ is constructed from $P'$ with the covering relation  defined as $v\lhd u$ for all $u\in M'$.
			\item  if $v$ is adjacent to the elements of $S_0'$, where $M'\subset S_0'\subset S'$, then $P$ is constructed with the covering relation defined as $v\lhd u'$ for all $u'\in S_0'\setminus M'$ and $v\lhd u\lhd u''$ for all $u''\in S'\setminus S_0'$ and some $u\in M'$.
			\item If $v$ is adjacent to all elements of $S'$, then $P$ is constructed from $P'$ with the covering relation $v\lhd u$ for all $u \in S'\setminus M'$.
		\end{itemize}
		It follows from construction that $P$ is a well-defined poset and that $G$ is isomorphic to $G_P$. That is, $G$ is a chordal C-I graph. Hence the lemma.
	\end{proof}

	A \emph{perfect elimination ordering (PEO)} is an ordering $\pi= v_1,\ldots,v_n$ of vertices in $G$  such that the neighborhood $N[v_i]$ of $v_i$ is a clique of the subgraph $G_{\{v_i,\ldots, v_n\} }$ induced by the vertices $\{v_i,\ldots, v_n\}$ of $G$.  The following characterizations of chordal graphs
	are well known.
	\begin{theorem}\cite{fulkerson}
		A graph $G$ is chordal if and only if $G$
		has a perfect elimination ordering.
	\end{theorem}
	\begin{lemma}\cite{Dirac-1}\label{lemma-chordal-dirac}
		Every chordal graph $G$ has a simplicial vertex. If $G$ is not complete, then it has two non-adjacent simplicial vertices.
	\end{lemma}
	
	In the following, we prove an interesting characterization of C-I chordal graphs as precisely those chordal graphs having exactly two independent simplicial vertices. 
	
	\begin{theorem}\label{chordal_2-simplicial}
		Let $G$ be a chordal graph. Then $G$ is a C-I graph if and only if $G$ is a complete graph or $G$ has exactly two independent simplicial vertices.
	\end{theorem}
	\begin{proof}
		Suppose $G$ is a chordal and C-I graph, then $G$ is a complete graph, or $G$ contains exactly two independent simplicial vertices by Lemma~\ref{CI3notsimplicial} and \ref{lemma-chordal-dirac}. 
  
  Conversely, every complete graph is a C-I graph. It remains to prove that if $G$ is chordal and $G$ has exactly two independent simplicial vertices, then $G$ is a C-I graph.
		
		Let $G$ be the chordal graph with the vertex set $V(G)=\{v_1,v_2,v_3,\ldots,v_n\}$. Since $G$ is chordal, there is a perfect elimination ordering(PEO), let $v_n,v_{n-1},\ldots,v_2,v_1$ be a PEO in $G$. By the definition of PEO on a chordal graph, the neighborhood $N(v_{n-i})$ of $v_{n-i}$ is a clique in the subgraph $G_{\{v_{
  n-\{i+1\}},\ldots, v_1\} }$, for $i=0,1,\ldots n-1$ and also the subgraphs $G_{\{v_{n-(i+1)},\ldots, v_1\} }$ are chordal for $i=0,1,\ldots n-1$.

According to Lemma~\ref{G-v_simplicial}, if a chordal graph has exactly two independent simplicial vertices, then the removal of a simplicial vertex results in either a chordal graph with exactly two independent simplicial vertices or a complete graph. We eliminate the simplicial vertices until we get the smallest non-trivial chordal graph containing exactly two independent simplicial vertices, say $G_k$. The removal of a simplicial vertex from $G_k$ results in a complete graph. That is, $G_{k-1}$ is a complete graph with $G_{k-1}=G_k\setminus {v_k}$, where $v_k$ is a simplicial vertex in $G_k$. Clearly, $G_{k-1}$ is a C-I graph being a complete graph. Now, we add the simplicial vertex $v_k$ to $G_{k-1}$ and form the graph $G_K$ by adding the deleted edges back. We continue this process by adding the eliminated simplicial vertices $v_{k+1}, v_{k+2}, \ldots v_{n}$ and deleted edges in the reverse order in the PEO  to the graphs $G_k, G_{k+1}, \ldots G_{n-1}$, respectively obtaining finally the graph $G_n \cong G$.  
		
		%Now, we add the simplicial vertices to the graphs $G_{k-1}$, $G_{k}, \ldots, G_{n-1}$, in reverse order in the PEO; that is, we add vertices in the order $v_k,v_{k+1},\ldots,v_n$ and it follows from the construction that we always get back the graph $G$.  
 It may be noted that we add edges so that the resulting graph obtained in each stage contains exactly two independent simplicial vertices.	
		
In particular, the graph $G_k$ is obtained by adding the simplicial vertex $v_k$ to $G_{k-1}$. The vertex $v_k$ is adjacent to some vertices of $G_{k-1}$. Let $C_{k-1}\subset V(G_{k-1})$ and $C_{k-1}'= V(G_{k-1})\setminus C_{k-1}$ such that $v_k$ is adjacent to every element of $C_{k-1}$ and not adjacent to any element of $C_{k-1}'$. We form a poset, say $P_0$ consisting of elements, $C_{k-1}'$, $v_k$, and $C_{k-1}$.  Now fix $C_{k-1}'$ as the set of minimal elements, $v_k$ covering every element in $C_{k-1}'$ and all elements of $C_{k-1}$ covering $v_k$. It is clear that the C-I graph of $P_0$ is isomorphic to $G_k$. That is, we obtain that $G_k$ is a C-I graph. It is in-fact a Ptolemaic C-I graph of the completely ranked poset with rank 2, where $C_{k-1}'$ represents elements of rank 0, $v_k$ represents the element of rank 1 and  $C_{k-1}$ represents elements of rank 2. Thus $G_k$ is a chordal C-I graph with exactly two simplicial vertices.
		
Now the graph $G_{k+1}$ is obtained by adding the simplicial vertex $v_{k+1}$ and the deleted edges to $G_k$ and from Lemma~\ref{G_k_chordal}, we obtain that $G_{k+1}$ is a chordal C-I graph having exactly two independent simplical vertices. By by the inductive arguments, we have that $G_{k+2}, G_{k+3},\ldots, G_{n}$ are chordal C-I graphs with exactly two independent simplicial vertices. Finally, we get that the graph $G_n$ is isomorphic to $G$ and hence $G$ is a C-I graph. Hence the theorem.

%$G_{k+2}$ is obtained by adding the simplicial vertex $v_{k+2}$ to $G_{k+1}$ and continuing like this, we obtain $G_n$ by adding the simplicial vertex $v_n$ to $G_{n-1}$. By  Lemma~\ref{G_k_chordal} and Corollary~\ref{coro_G-v_2_indpendent}, we have that $
		
	\end{proof}

\subsection{Recognition of chordal C-I graphs}
\begin{figure}[htb]	
		\epsfxsize=10truecm \slika{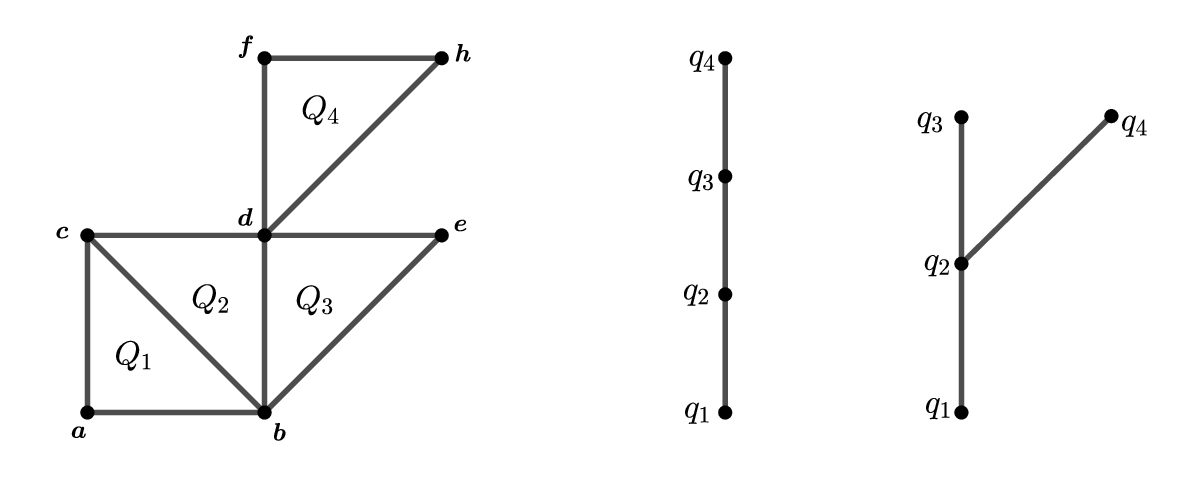}\vspace{-.6cm} \caption{A chordal graph and two of its clique trees}\label{fig:clique-tree}
	\end{figure}

	In this section, we design an algorithm for recognizing chordal C-I graphs. For this, we use the clique tree representation of a chordal graph. Since this algorithm computes a clique tree, if the graph is a chordal graph explicitly, it turns out that the algorithm provides a certification. 
	It is well known that a chordal graph $G$ can be characterized as the intersection graph of subtrees of some tree $T$. Such a tree is known as the clique tree of $G$. This well-known theorem is proved independently by Buneman \cite{Buneman} and Gavril \cite{gravril}.  A tree representation of a chordal graph $G$ is a pair $(T,\mathscr{F})$ where $T$ is a tree and $\mathscr{F}$ is a family of subtrees of $T$ such that the intersection graph of $\mathscr{F}$ is isomorphic to $G$.  Further Gavril \cite{gravril} has
	shown that given a chordal graph $G$, it is possible to construct a tree $T$ with vertex set $K=\{q_1, q_2,\ldots,q_r\}$ where $q_i$ corresponds to the maximal clique $Q_i$ of $G$, such that
	$(T, \{R_{v_1}, R_{v_2},\ldots,R_{v_n}\})$ is a tree representation of $G$. Here, each $R_{v_i}$, $1\leq i\leq n$, is the set of maximal cliques that contain the vertex $v_i$, that is, $R_{v_i} = \{q_j \mid v_i\in Q_j\}$. Such a tree
	representation of $G$ is called a \emph{clique tree} of $G$. The clique tree $T$ of $G$ need not be unique, for the sets $R_i$ determine the edges of $T$, not necessarily in a unique manner. %An undirected graph is chordal if and only if it has a clique tree. 
	Fig.\ref{fig:clique-tree} shows a chordal graph and two of its clique trees.
	
	\begin{theorem}[Gavril \cite{gravril} and Buneman \cite{Buneman}]
		The following propositions are equivalent:
		\begin{enumerate}
			\item $G$ is a chordal graph.
			\item  $G$ is the intersection graph of a family of subtrees of a tree.
			\item There exists a tree $T$(called clique tree) with vertex set $\{q_1, q_2,\ldots,q_r\}$ such that for each vertex $v$,
			the set $R_v = \{q_i\mid v\in Q_i\}$ induces a subtree of $T$. Here $q_i$ represents the maximal clique $Q_i$.
		\end{enumerate}
	\end{theorem}

	From the above description, we have that a graph $G$ is a chordal graph if and only if there exists a clique tree $T$ with maximal cliques of $G$ as the vertices of $T$ and any two cliques containing $v\in V(G)$ are either adjacent in $T$ or connected by a path of cliques that contain $v$. The Lemma below is also very useful in understanding the structure of chordal graphs and their simplicial vertices.

	\begin{lemma}\cite{Blair}\label{clique-tree-simplicial}
		A vertex is simplicial if and only if it belongs to precisely one maximal clique.
	\end{lemma}

	Even though the result of the following lemma is a known property of chordal graphs, we require this for establishing the recognition of chordal graphs which are C-I graphs, and hence we state it below. 
	\begin{lemma}\label{leaf-simplicial}
		Every leaf node of a clique tree of a chordal graph contains a simplicial vertex.
	\end{lemma}
	\begin{proof}
		Let $T$ be a clique tree of a chordal graph $G$. Let $C_k$ be a leaf node of $T$. Since $C_k$ is a leaf node of $T$ there is a node $C_{k'}$ in $T$ such that $C_kC_{k'}\in E(T)$. Since $T$ is a clique tree of $G$ every node of $T$ is a maximal clique of $G$, and hence there is a vertex $v$ of $G$ in $C_k$ but not in $C_{k'}$. From the definition of the clique tree, it is clear that  $v$ is not in any other nodes of $T$ (cliques). Then it follows that $v$ is a simplicial vertex by Lemma~\ref{clique-tree-simplicial}.
	\end{proof}
	
 The following lemma will also be useful for the recognition of chordal C-I graphs. 
 
	\begin{lemma}\label{path-simplicial}
		Let $G$ be a chordal graph with the clique tree $T$ being a path. Let $C_i$ be an internal node of $T$, $C_{i+1}$, the parent node of $C_i$, and $C_{i-1}$ the child node of $C_i$ in $T$. If $v$ is vertex of $G$ such that $v\in C_i\setminus C_{i+1}$ and $v\in C_i\setminus C_{i-1}$, then $v$ is a simplicial vertex in $G$. 
	\end{lemma}
	\begin{proof}
		Since  $v\in C_i\setminus C_{i+1}$ and $v\in C_i\setminus C_{i-1}$ and $T$ is a path, it follows from the definition of a clique tree that $v$ is not in any other cliques in $T$. Then by Lemma~\ref{clique-tree-simplicial}, $v$ is a simplicial vertex.
	\end{proof}

	Based on Theorem~ \ref{chordal_2-simplicial} and Lemma~\ref{leaf-simplicial} and \ref{path-simplicial}, we formulate the following algorithm for recognizing a chordal C-I graph $G$. The correctness of the algorithm follows from these results.

	\begin{small}
		\begin{algorithm}[H] \label{algo-chordal}
			\DontPrintSemicolon
			\caption{  Algorithm for recognizing given graph $G$ is chordal C-I graph $G$ or not.}
			\SetAlgoLined
			%\KwResult{Write here the result}
			\SetKwInOut{Input}{Input}\SetKwInOut{Output}{Output}
			\vspace{.25cm}
			\Input{$G$ be a connected graph with $|V(G)|=n$ and $|E(G)|=m$}
			\Output{$G$ is  chordal C-I graph or not } 
			\vspace{.25cm}
			\begin{itemize}
				%\item[Step 0: ]  If $G$ is a complete graph then return $G$ chordal C-I graph and Stop. Otherwise, go to Step 1.
				\item[Step 1: ] 
				Apply the perfect elimination ordering (PEO) on $G$ and check whether $G$ is a chordal graph or not.\\ If there is no PEO then return $G$ not chordal and stop. Otherwise, go to Step 2.
				\item[Step 2:] Using the PEO find the clique tree $T$ of the given graph $G$ using Blair-Peyton algorithm.\\ If the clique tree has more than 2 leaf nodes, then return $G$, not chordal C-I graph, and stop. \\ Otherwise, go to Step 3. (Now the clique tree is a path, and the vertices are in order $C_1, C_2,\ldots, C_k$)
				
				\item[Step 3:] Check whether there is an element in $C_i$ that is not in $C_{i+1}$ and $C_{i-1}$ for $i=2,3,\ldots,k-1$.\\ If such an element exists at any stage, then stop and return $G$ as not a chordal C-I graph.\\ Otherwise, return $G$  as a chordal C-I graph.

			\end{itemize}
		\end{algorithm}
	\end{small}
	
	The time complexity of Algorithm \ref{algo-chordal} can be analyzed as follows. 

Using Lex-BFS or MCS, a PEO can be found in  $O(n + m)$ time. In Step 2, the Blair-Peyton algorithm is used to construct a clique tree from a chordal graph based on the PEO \cite{Blair}. This algorithm also runs in linear time, $O(n + m)$. Additionally, checking whether there are more than two leaf nodes can be done in constant time. 	
 
In Step 3, let $k$ be the length of the clique tree which is a path $P$. Then, in the worst case, the size of the sets $C_i$ is of the order $O(n/k)$. To check whether an element $v\in C_i\setminus C_{i+1}$ and $v\in C_i\setminus C_{i-1}$ takes $O(3n/k)$ time by using a hash set operation. Since there are at most $k-2$ internal nodes, the total time for Step 3 is $O(n)$.
	Now, the complexity of Algorithm \ref{algo-chordal} is $O(n+m)$.

 According to Theorem~\ref{chordal_2-simplicial}, a graph $G$ is a chordal C-I graph if and only if $G$ is either a complete graph or has exactly two independent simplicial vertices. Lemma~\ref{clique-tree-simplicial} and Lemma~\ref{leaf-simplicial} further imply that the clique tree of a chordal C-I graph contains exactly two leaf nodes, meaning that the clique tree forms a path. The next step is to certify that none of the internal vertices of the clique tree contain simplicial vertices. Lemma~\ref{clique-tree-simplicial} and Lemma~\ref{path-simplicial} address this in Step 3 of the algorithm.
 
 Since the Blair-Peyton algorithm explicitly computes the clique tree when $G$ is a chordal graph and after Steps 2 and 3 of the algorithm outputs a clique tree which is a path and contains no internal vertices with simplicial vertices, if the given graph is a chordal C-I graph. So the algorithm actually returns a chordal C-I graph if $G$ is such a graph, which is a certification model.
 
 %The Blair-Peyton algorithm constructs a clique tree from a chordal graph, which serves as a certificate of chordality \cite{Blair}. First, determine whether the graph is chordal using either Lex-BFS or MCS. Then, apply the Blair-Peyton algorithm to construct the clique tree. This clique tree acts as a certificate because it verifies the structure of the chordal graph.

%According to Theorem~\ref{chordal_2-simplicial}, a graph $G$ is a chordal C-I graph if and only if $G$ is either a complete graph or has exactly two independent simplicial vertices. Lemma~\ref{clique-tree-simplicial} and Lemma~\ref{leaf-simplicial} further imply that the clique tree of a chordal C-I graph contains exactly two leaf nodes, meaning that the clique tree forms a path. The next step is to certify that none of the internal vertices of the clique tree contain simplicial vertices. Lemma~\ref{leaf-simplicial} and Lemma~\ref{path-simplicial} address this in Step 3. This process certifies chordal C-I graphs and runs in linear time.

	\section{Structure of C-I cographs}\label{Section_cograph}
	
	In this section, we present an algorithm for recognizing C-I cograph.
	\begin{theorem}\cite{cograph}
		Let $G$ be a chordal cograph. Then $G$ is a C-I graph if and only if
		$G$ is a connected graph that contains at most two maximal cliques.
		\label{chordal cograph}
	\end{theorem}
	
\begin{figure}[htb]	
		\epsfxsize=5truecm \slika{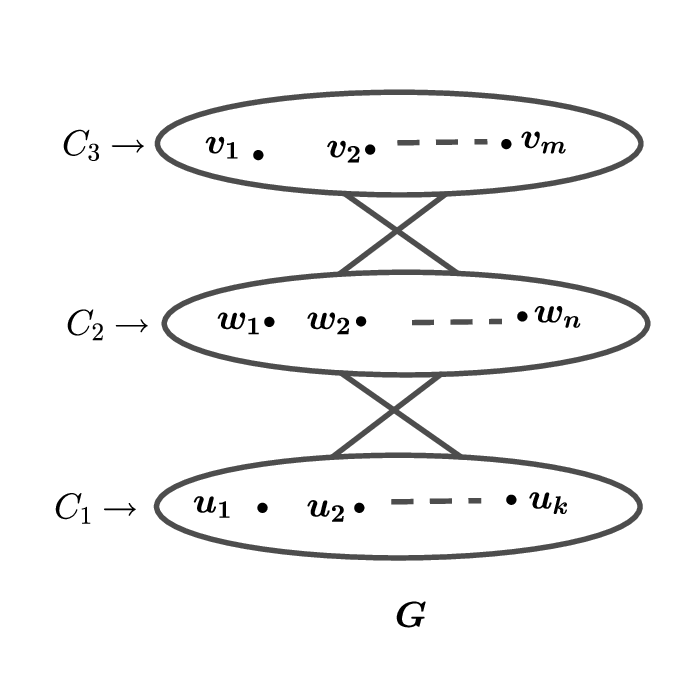}\vspace{-.75cm} \caption{General structure of a chordal C-I cograph $G$}\label{fig:Chordal_cograph_structure}
	\end{figure}

	From Theorem \ref{chordal cograph}, we get the following remark
	\begin{remark}\label{remark_chordal_cograph}
		A graph $G$ is a chordal C-I cograph if and only if there exists three pairwise disjoint sets $C_1, C_2$ and $C_3$ such that $V(G)=C_1\cup C_2\cup C_3$, and $x,y \in V(G)$ are adjacent in $G$ if and only if $x,y\in C_1\cup C_2$ or $x,y\in C_2\cup C_3$.
	\end{remark} 
	That is, $C_1, C_2, C_3, C_1\cup C_2$ and $C_2\cup C_3$ form cliques and the graph has no other edges. So in the C-I chordal cograph, if $xy,xz\notin E(G)$ then $yz\in E(G)$. It is clear that a chordal C-I cograph is a Ptolemaic C-I graph as it is a C-I graph of a completely ranked poset of height 3. 
	The structure of an arbitrary chordal C-I cograph is shown in Fig.~\ref{fig:Chordal_cograph_structure}.
	\begin{theorem}\cite{cograph}
		A graph $G$ is a C-I cograph if and only if $G$ is the join of chordal C-I cographs.
		\label{thm_cograph}
	\end{theorem}
	From Theorem~\ref{thm_cograph}, we get the following remark.
	\begin{remark}\label{remark_cograph}
		Let $G$ be a C-I cograph such that $G=G_1\vee G_2 \vee\cdots\vee G_k$, where $G_i$'s are chordal C-I cographs. By Remark~\ref{remark_chordal_cograph}, there exists a pairwise disjoint set $C_1^{i}, C_2^{i}$ and $C_3^{i}$ such that $V(G_i)=C_1^{i}\cup C_2^{i}\cup C_3^{i}$ and $xy\in E(G_i)$ if and only if $x,y\in C_1^{i}\cup C_2^{i}$ or $x,y\in C_2^{i}\cup C_3^{i}$. Thus in $G$, $\bigcup\limits_{i=1}^{k} C_2^{i}$ are universal vertices and for $u,v\in V(G)$ such that $uv\notin E(G)$ if and only if both $u$ and $v$ in $V(G_i)$  with $u\in C_1^{i}$ and $v\in C_3^{i}$, for some $i$. Thus, for any $u,v,w\in V(G)$ with $uv,uw\notin E(G)$, then $vw\in E(G)$. 
	\end{remark}

\begin{figure}[htb]	
		\epsfxsize=14.5truecm \slika{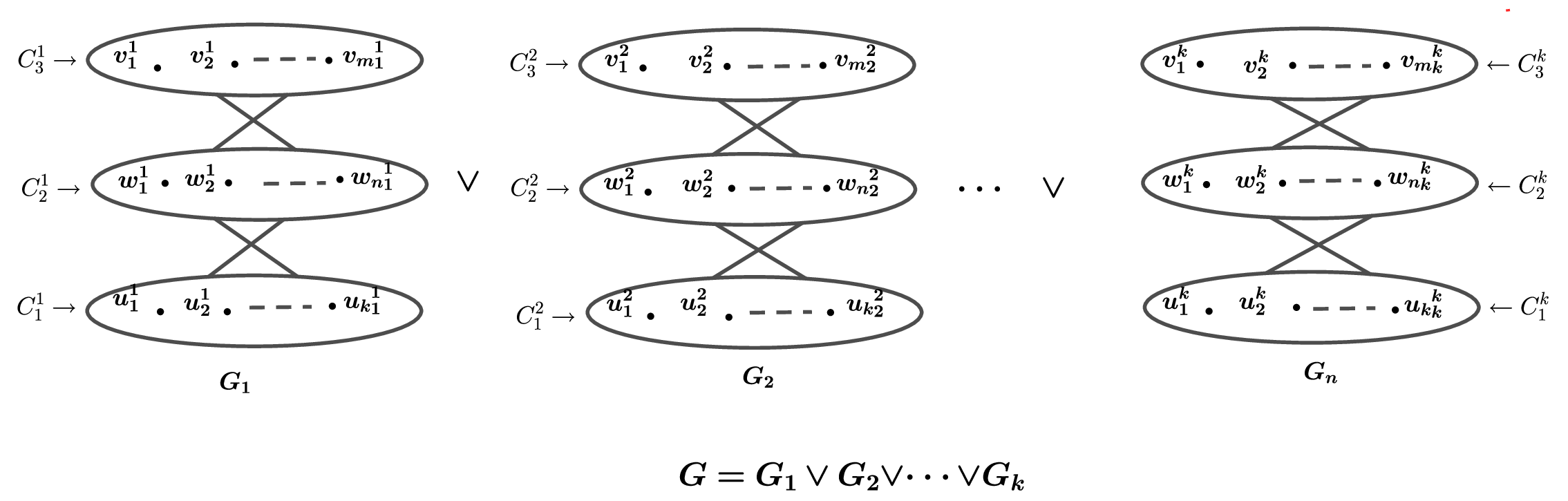}\vspace{-.4cm} \caption{General structure of a C-I cograph $G$}\label{fig:cograph_structure}
	\end{figure}

	The structure of an arbitrary C-I cograph is depicted in Fig.~\ref{fig:cograph_structure}.
	\begin{lemma}\label{algo_cograph_poset} 
		Let $G$ be a C-I cograph then there exists a poset  $P=P_1+P_2 + \cdots + P_r$, where $P_i$'s are completely ranked poset of rank 2 such that $G\cong G_P$.
	\end{lemma}
	\begin{proof}
		Let $G$ be a C-I cograph. Then $G$ is the join of chordal C-I cographs. That is, let $G=G_1\vee G_2\vee \dots G_r$. Since $G_i$'s are chordal cographs, $G_i$'s are C-I graphs of the completely ranked poset $P_i$'s of rank at most 2. Hence $P=P_1+ P_2+\cdots + P_r$.
	\end{proof}

	\begin{figure}[htb]	
		\epsfxsize=4truecm \slika{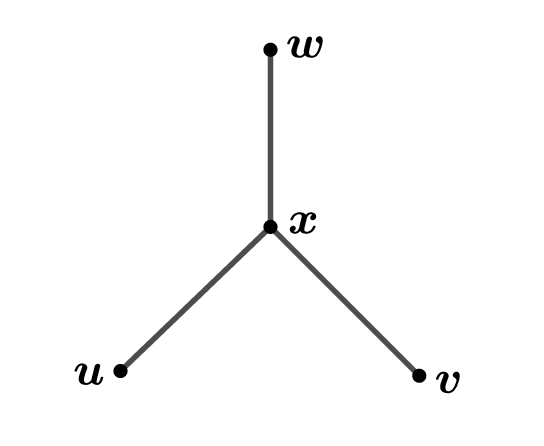}\vspace{-.65cm} \caption{Claw }\label{fig:claw}
	\end{figure}
\begin{figure}[htb]	
		\epsfxsize=14.5truecm \slika{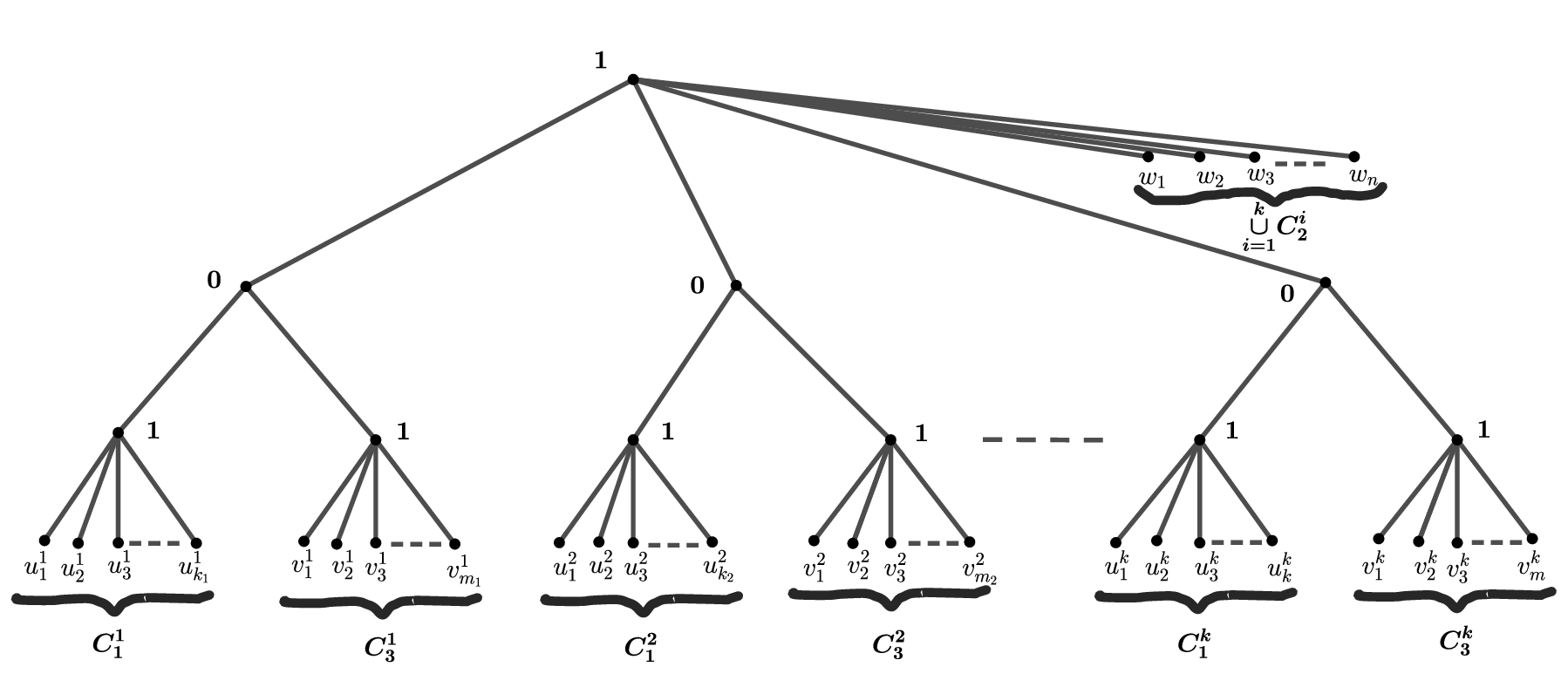}\vspace{-.5cm} \caption{Cotree for C-I cograph $(\mathcal{T_C})$ }\label{fig:cograph_cotree}
	\end{figure}

	\begin{lemma}\label{lemma:claw}
		If $G$ is a C-I cograph then $G$ is a claw-free graph. 
	\end{lemma}
	\begin{proof}
		Let $G$ be a C-I cograph and $G=G_1\vee G_2\vee\cdots \vee G_k$, where $G_i$'s are chordal C-I cograph for $i=1,2,\ldots k$.
		Suppose that $G$ contains an induced claw by the vertices $\{u,v,w,x\}\subseteq V(G)$ with edges $ux,vx,wx\in E(G)$ and $uv,vw,uw\notin E(G)$ (see Fig.\ref{fig:cograph_cotree}). Since $G=G_1\vee G_2\vee\cdots \vee G_k$, the non-adjacent vertices lie in the same $G_j$ for some $j\in \{1,2,\ldots ,k\}$. That is, $u,v,w \in V(G_j)$. Since $G_j$ is a chordal C-I cograph, if $uv,uw\notin E(G_j)$ then $vw\in E(G_j)$, which is a contradiction. Hence, the C-I cograph is a claw-free graph.
	\end{proof}

	A \emph{cotree} is a tree in which the internal nodes are labelled with the numbers 0 and 1. Every cotree $T$ defines a cograph $G$ having the leaves of $T$ as vertices, and in which the subtree rooted at each node of $T$ corresponds to the induced subgraph in $G$ defined by the set of leaves descending from that node.  
	A subtree rooted at a node labelled 0 corresponds to the union of the subgraphs defined by the children of that node and a node labelled 1 corresponds to the join of the subgraphs defined by the children of that node.
	The cotree satisfies the property that, on every root-to-leaf path, leaves are the vertices of the graph, the labels of the internal nodes alternate between 0 and 1, and every internal node has at least two children. The cotree can be easily obtained from any tree labelled with such $0\slash1$ $T$ by coalescing all pairs of child-parent nodes in $T$ having the same label or where the parent has only one child. Vertices $x$ and $y$ of $G$ are adjacent in $G$ if and only if their lowest common ancestor(LCA) in the cotree is labelled 1. This representation is unique and every cograph can be represented in this way by a cotree \cite{corneil_D_G}. Since $G$ is a connected graph, the root of the cotree is a 1-node. 
	
	It is known that cographs have a unique tree representation, called a cotree. Using the cotree, it is possible to design very fast polynomial-time algorithms for problems that are intractable for graphs in general. Such problems include chromatic number, clique determination, clustering, minimum weight domination, isomorphism, minimum fill-in, and Hamiltonicity. A linear-time cograph recognition algorithm, such as those in \cite{habib-2008,corneil,habib-2005}, which also builds a cotree, a data structure that fully encodes the cograph.
	
\begin{figure}[htb]	
		\epsfxsize=14.5truecm \slika{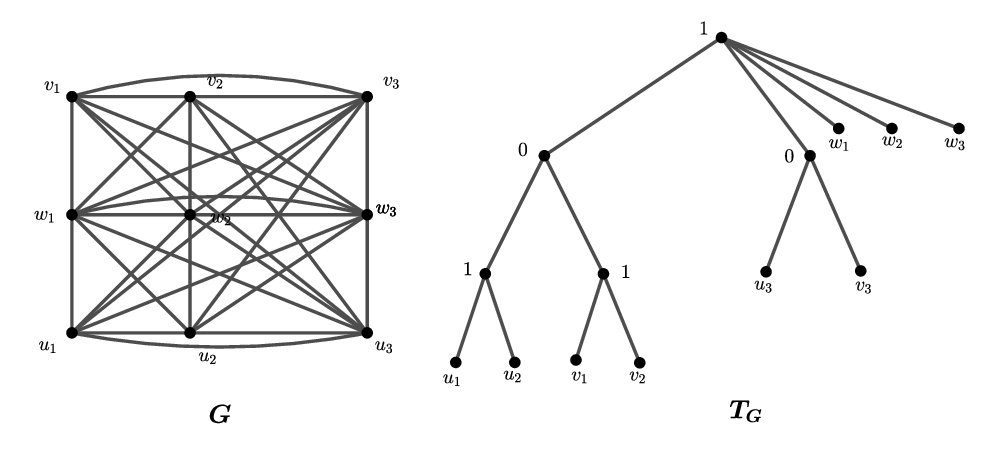}\vspace{-.75cm} \caption{$G$ be a cograph and $T_G$ be the cotree of $G$}\label{fig:cograph_example}
	\end{figure}

	The Fig.~\ref{fig:cograph_example}  illustrates a cograph $G$  and an embedding of the corresponding cotree $T_G$. The leaves of $T_G$ represent the set of vertex $V(G)$, and each internal node signifies the union (0) or the joining (1) operations in the children. The significance of the 0(1) nodes is captured by the fact that $xy\in E(G)$  if and only if $LCA(x, y)$ is a 1 node, as shown in Fig.~\ref{fig:cograph_example}.
	
	The general structure of the cotree of C-I cographs is shown in Fig.~\ref{fig:cograph_cotree}, and we denote the family of cotree of C-I cograph as $\mathcal{T_C}$. The children of the root node are referred to as \emph{level 2} nodes, while the children of level 2 nodes are termed \emph{level 3} nodes.
	Let $G$ be a C-I cograph, $G=G_1\vee G_2 \vee\cdots\vee G_k$, where $G_i$'s are chordal C-I cographs. Then by Remark~\ref{remark_cograph},
	two vertices $x$ and $y$ are not adjacent in $G$ if and only if $x\in C_1^i$ and $y \in C_3^i$ for $i=1,2,\ldots, k$. Therefore, in the cotree, the lowest common ancestor(LCA) of $x$ and $y$ is 0-node. That is, $LCA(x,y)=0$ for $x\in C_1^i$ and $y \in C_3^i$. In other cases, LCA is a 1-node. That is, $LCA(x,y)=0$ for $x\notin C_1^i$ and $y \notin C_3^i$.

	\begin{lemma}\label{lemma_cotree_condition}
		The cotree of a C-I cograph which is not a complete graph satisfies the following properties: 
		\begin{itemize}
			\item[(i)] The root node (which is a 1-node) of $\mathcal{T_C}$ has leaf nodes and $0$-nodes as children. the number of 0-nodes (which are children of the root node) is always less than or equal to the number of leaf nodes (which are children of the root node).%the number of $0$-node children is less than the number of leaf node children.
			\item[(ii)] Every 0-node has exactly two children (both the children are 1-nodes or both the children are leaves nodes or one child is a leaf node and another child is 1-node).
			\item[(iii)]The $1$-nodes other than the root node have only leaf nodes as children.
			\item[(iv)] The cotree of $G$ belongs to the family of $\mathcal{T_C}$.
			
		\end{itemize}

	\end{lemma}
	\begin{proof}
		Let $G$ be a C-I cograph which is not a complete graph.
		\begin{itemize}
			\item[(i)] Since $G$ is not a complete C-I graph, it is either a chordal cograph or joins of chordal cographs. Since there are non-adjacent vertices, there exist some 0-nodes as children of the root node. From the structure of the C-I cograph, there exist universal vertices, which should be leaf nodes of the root node. From Theorem~\ref{thm_cograph}, $G=G_1\vee G_2\vee\cdots \vee G_k$, where $G_i$ are chordal C-I cographs. The leaf nodes of the subtree rooted by each 0-node contain the non-universal vertices of each $G_i$. That is, if there are $k$ 0-nodes, then since the universal vertices of each $G_i$, are also universal vertices of $G$, and since each $G_i$ has at least one universal vertex, the number of 0-nodes of the root node must be at most the number of leaf nodes as a child nodes. Hence (i) follows.
			\item[(ii)] Suppose there are $0$-nodes having more than two children. Let $0'$ be a $0$-node having three child nodes. Let $u,v,w$ be the three leaf nodes from the different branches of the subtree with the root node as the $0'$-node. Then by the structure of the C-I cograph, we have that $uv,vw,uw \notin E(G)$. Let $x$ be a leaf node of the root of the cotree $\mathcal{T_C}$. Then  $x$ is adjacent to $u,v$ and $w$. That is,  $\{u,v,w,x\}$ induce a claw in $G$ which is a contradiction to Lemma~\ref{lemma:claw}. So every $0$-node has no more than two children. In a cotree, every node has a minimum of two child nodes. Hence every $0$-node has exactly two child nodes as 1-nodes. Hence (ii) follows.
			\item[(iii)]Suppose that there is a $1$-node other than the root node which has a 0-node, say $0'$. Consider a subtree $T$ of $\mathcal{T_C}$ with rooted in $0'$.  Let $u$ and $v$ be the vertices of $G$ such that the $LCA (u,v)$ is $0'$. The vertices $u$ and $v$ exist since the subtree $T$ has exactly two branches. Now consider the subtree $T'$ of $\mathcal{T_C}$ containing $T$ with root as a 0-node, say $0''$ different from $0'$. Let $w$ be a vertex of $G$ such that the $LCA(u,v,w)$ is $0''$ ( the vertex $w$ exists since there are two children from $0''$ of $T'$). Clearly, the vertices $u,v, w$ are mutually non-adjacent. Now any vertex $x$ of $G$ which are leaves of the root node 1 of $\mathcal{T_C}$ is adjacent to all of $u,v,w$ so that $u,v,w,x$ form an induced claw, a contradiction. Hence (iii) follows.
			\item[(iv)] Follows from (i),(ii) and (iii).
		\end{itemize}
	\end{proof}

	\begin{theorem}\label{thm-cotree}
		A cograph $G$ is a C-I cograph if and only if the cotree of $G$ belongs to the family $\mathcal{T_C}$.  
	\end{theorem}
	\begin{proof}
		It is clear that when $G$ is a complete graph if and only if the cotree is the tree with every vertex of $G$ as a child of the root node and which is a part of $\mathcal{T_C}$.
		
		Let $G$ be a C-I cograph which is not a complete graph. Then by Lemma~ \ref{lemma_cotree_condition}, cotree of $G$ is of the form $\mathcal{T_C}$
		
		Conversely, we need to prove that if $\mathcal{T_C}$ is the cotree of a cograph $G$ then $G$ is a C-I graph. Consider the set of leaves of the cotree $\mathcal{T_C}$ ( Refer Figure~\ref{fig:cograph_cotree} ). It follows by the definition of cotree of a cograph $G$ that the sets $C_i^j$'s are cliques of $G$, for $i=1,3 $ and $j=1,2,\ldots,k$ and $\bigcup\limits_{i=1}^{k} C_2^{i}=\{w_1,w_2,\ldots,w_m\}$ are universal vertices. Now $C_1^j$ is adjacent to every vertex in $G$ except $C_3^j$ and $C_3^j$ is adjacent to every vertex in $G$ except $C_1^j$ for $j=1,2,\ldots,k$. Then by Remark~\ref{remark_cograph}, $G$ is a C-I cograph. 
		
	\end{proof}
	
	In \cite{corneil}, Corneil et al. presented a linear-time algorithm for recognizing cographs and constructing their cotree representation. Using that algorithm, we get the cotree.
	\begin{theorem}\cite{corneil}\label{thm_corneil}
		The family of cographs possesses a linear-time recognition algorithm by the construction of the cotrees in linear-time.
	\end{theorem}
	\begin{small}
	\begin{algorithm}[H] \label{algorithm-cograph}
		\DontPrintSemicolon
		\caption{ Algorithm for recognizing C-I cograph}
		\SetAlgoLined
		%\KwResult{Write here the result}
		\SetKwInOut{Input}{Input}\SetKwInOut{Output}{Output}
		\vspace{.25cm}
		\Input{A connected graph $G$ with $|V(G)|=n$}
		%\Input{The cotree $T$ for the cograph $G$ using the recognition algorithm in Theorem~\cite{} **.}
		\Output{$G$ is a C-I cograph or not. If G is a C-I cograph then find a poset $P$
such that $G_P\cong G$.}
\vspace{.25cm}
		\begin{itemize}
			\item[Step 0:   ] Using the recognition algorithm  in Theorem~\ref{thm_corneil}  determine the cotree if $G$ is a cograph and go to Step 1. Otherwise, stop and return $G$ is not cograph
			\item[Step 1:   ] Perform BFS from the root node 1 of the cotree.		
			\item[Step 2:   ] If BFS completes in one step by reaching all the leaf nodes (if all the adjacent nodes of the root node are\\ leaf nodes), then the graph is a complete graph and it is C-I cograph; stop, and we are done.\\ Otherwise, go to Step 3.
  			\item[Step 3(i): ] After Step 2, only 0-nodes are obtained. Then $G$ is not a C-I cograph and Stop.
     \item[Step 3(ii):] After Step 2, both 0-nodes and leaf nodes are obtained. Check whether the number of zero nodes is\\ greater than the number of leaf nodes. If yes, then $G$ is not a C-I cograph and stop.\\ Otherwise goto Step 4
     
   %  \textcolor{red}{( that is, in Step 2, we have two possibilities (i) Both 0-nodes and leaf nodes are obtained (ii) only 0-nodes are obtained.) If all the neighbours of the root nodes are only $0$-nodes, then stop, since $G$ is not a C-I graph.  or the number of 0-nodes is greater than the number of leaf nodes of the root node, then $G$ is not a C-I cograph and stop. Otherwise, go to Step 4.}
			
			\item[Step 4:   ] Continue the searching at level 2.  If $0$-nodes has more than two children as $1$-nodes then $G$ is not a C-I cograph and stop. Otherwise, go to Step 5.
			\item[Step 5:   ] Continue searching at level 3. In level 3, if any $1$-node has $0$-node as a child, then $G$ is not a C-I cograph\\ and stop. If all the children of the $1$-nodes are the leaf node then we are done. The resulting graph is a C-I cograph.

     \item[Step 6:   ] Label the leaf nodes of each 0-node ~branch in level 2 as $C_1^i$ and $C_3^i$. Partition the leaf nodes from\\ the rooted node 1 into the number of 0-nodes. That is, the leaf nodes are partitioned into $r$ sets,  $C_2^i$ for $i=1,2,\ldots,r$, where $r$ denotes the number of 0-nodes at level 2. 

     \item[Step 7:   ] Construct  completely ranked posets $P_i$ with rank 0 elements as $C_1^i$, rank 1 elements as $C_2^i$, and rank 2 elements as $C_3^i$ for $i=1,2,3,\ldots,r$.\\ Return  $P=P_1+P_2+\cdots+P_r$ as one of the poset of the C-I cograph $G$.
		\end{itemize}
	\end{algorithm}
\end{small}

%After Step 5, label the leaf nodes of each 0-node in level 2. Let $C_1^i$ and $C_3^i$ be the leaf nodes of the $i^{th}$ 0-node (there is a possibility that $C_3^i$ is empty). Partition the leaf nodes within the rooted node into $C_2^i$ for $i=1,2,\ldots,r$, where $r$ denotes the number of 0-nodes at level 2. Then, construct completely ranked posets $P_i$ with rank 0 elements as $C_1^i$, rank 1 elements as $C_2^i$, and rank 2 elements as $C_3^i$ for $i=1,2,3,\ldots,r$. Then, $P=P_1+P_2+\cdots+P_r$ is a poset of the C-I cograph (by Lemma \ref{algo_cograph_poset} and Theorem \ref{thm-cotree}.

\textbf{Proof of correctness}: 	The correctness of the Algorithm~\ref{algorithm-cograph} up to step 5 follows from the structure of cotree stated in Lemma \ref{lemma_cotree_condition} and Theorem \ref{thm-cotree}, and the related facts mentioned above. If after step 5 results in a C-I cograph, then the poset $P=P_1+P_2+\cdots+P_r$ is certification and it follows from Remark \ref{remark_cograph}, Lemma \ref{algo_cograph_poset} and Theorem \ref{thm-cotree} that it is one of the completely ranked poset corresponding to the input graph $G$.

	Now the size of the cotree of a C-I cograph can be estimated as follows. Let $G$ be a C-I cograph with $n$ vertices and $G=G_1\vee G_2\vee\cdots G_k$. Then from the structure of the cotree $T$ of $G$, $T$ contains leaf nodes (which are the $n$-vertices of $G$), 0-nodes, and 1-nodes. Then $T$ contains $k$ 0-nodes, and each 0-node has exactly 2 children. Therefore $T$ has at most  $1+k+2k+n$ vertices. The worst case for $k$ is when $k= \frac{n}{3}$. Therefore in the worst case, the total number of vertices in $T$ is $1+\frac{n}{3}+\frac{2n}{3}+3$, which is of the order of $O(n)$.

	The time complexity of the recognition algorithm for cographs and finding its cotree by Theorem~\ref{thm_corneil} is linear. Now the complexity of Algorithm~\ref{algorithm-cograph} is linear in the size of the cotree, which is $O(n)$ in the worst case since the complexity of the BFS can be performed on the cotree in linear-time. 
	
	\subsubsection*{Acknowledgment:}Arun Anil  acknowledges the financial support from the University of Kerala, for providing University Post Doctoral Fellowship (Ac EVII 5576/2023/UOK dated 30/06/2023).


\begin{thebibliography}{99}
		
			
		\bibitem{am-03}Anil,~A., Changat,~M.:  Comparability Graphs Among Cover-Incomparability Graphs. In: Balachandran, N., Inkulu, R. (Eds.), CALDAM 2022, \textit{Algorithms and Discrete Applied Mathematics},  LNCS 13179, pp. 48--61, Springer, Cham (2022). https://doi.org/10.1007/978-3-030-95018-7\_5.
		
		\bibitem{am-04} Anil,~A., Changat,~M.: Composition and Product of Cover-Incomparability Graphs. \textit{The Art of Discrete and Applied Mathematics} 6, \#P2.09(2023). https://doi.org/10.26493/2590-9770.1515.af9
	
		
		\bibitem{am-02} Anil,~A., Changat,~M.: Ptolemaic and Chordal Cover-Incomparability Graphs. \textit{Order} 39, 29--43(2022). https://doi.org/10.1007/s11083-021-09551-w
		
		\bibitem{amtb-01} Anil,~A., Changat,~M., Gologranc,~T., Sukumaran,~B.: Ptolemaic and planar cover-incomparability graphs. \textit{Order} 38, 421--439(2021). https://doi.org/10.1007/s11083-021-09549-4
		
		
		\bibitem{Blair} Blair, J.R.S., Peyton, B.:  An Introduction to Chordal Graphs and Clique Trees. In: George, A., Gilbert, J.R., Liu, J.W.H. (Eds.), \textit{Graph Theory and Sparse Matrix Computation}. The IMA Volumes in Mathematics and its Applications, vol 56, pp. 1-29 Springer(1993). %https://doi.org/10.1007/978-1-4613-8369-7_1
  \bibitem{janbok-02} Bok, J., Maxov\'{a}, J.: Characterizing Subclasses of Cover-Incomparability Graphs by Forbidden Subposets. \textit{Order} 36, 349--358(2019).

  \bibitem{bcgmm-09} Bre\v{s}ar,~B.,  Changat,~M., Gologranc,~T., Mathew,~J.,  Mathews,~A.: Cover-incomparability graphs and chordal graphs. \textit{Discrete Appl. Math.} 158,  1752 --1759(2010).
		
		\bibitem{bckkmm-07} Bre\v{s}ar,~B.,  Changat,~M., Klav\v{z}ar,~S., Kov\v{s}e,~M.,  Mathew,~J., Mathews,~A.: Cover-incomparability graphs of posets. \textit{Order} 25, 335--347(2008).

  
		
		

  \bibitem{cograph} Bre\v{s}ar,~B.,  Gologranc,~T., Changat,~M., Sukumaran,~B.: Cographs Which are Cover-Incomparability Graphs of Posets. \textit{Order} 32, 179--187(2015).
		
		
		\bibitem{habib-2008} Bretscher, A., Corneil, D.G., Habib, M., Paul,  C.: A simple linear time LexBFS cograph recognition algorithm. \textit{SIAM J. Discrete Math.} 22 1277--1296(2008).
		
		
		
		\bibitem{bw-00} Brightwell,~G., West,~W.B.: Partially ordered sets. Chapter 11 in \textit{Handbook of Discrete and Combinatorial Mathematics} (K.H.~Rosen, ed.) CRC Press, Boca Raton, 717--752(2000).
		
		\bibitem{Buneman} Buneman, P.: A characterization of rigid circuit graphs. \textit{Discrete Math.} 9, 205--212(1974).
		
		\bibitem{corneil_D_G} Corneil, D.G., Lerchs, H.,  Stewart Burlingham, L.: Complement reducible graphs. \textit{Discrete Appl. Math.} 3, 163--174(1981). 
		
		\bibitem{corneil} Corneil, D.G., Perl, Y., Stewart, L.K.: A linear recognition algorithm for cographs.  \textit{SIAM J. Computing} 14, 926--934(1985).
		
		\bibitem{habib-c-1} Cournier, A., Habib, M.: A new linear algorithm for modular decomposition. In: Tison,S. (Ed.), \textit{$19^{th}$ International Colloquium Trees in Algebra and Programming}, CAAP'94, Vol. 787, LNCS, Springer, Berlin, 68--82(1994).
		
		\bibitem{Dahlhaus-2001} Dahlhaus, E., Gustedt, J., McConnell, R.M.: Efficient and practical algorithms for sequential modular decomposition. \textit{J. Algorithms} 41(2), 360--387(2001).
		
		\bibitem{Dirac-1} Dirac, G., A.: On rigid circuit graphs. \textit{Abh. Math. Sem. Univ. HamDurg}, 25, 71-76(1961).
		
		\bibitem{fulkerson} Fulkerson, D.R., Gross, O.A.: Incidence matrices and interval graphs. \textit{Pacific
		J. Math.} 15, 835-855(1965).
		
		\bibitem{Gallai} Gallai,~T.: Transitiv orientierbare Graphen. \textit{Acta Math. Acad. Sci. Hung.} 18, 25--66(1967).
		
		\bibitem{gravril} Gavril,~F.:  The intersection graphs of a path in a tree are exactly the chordal graphs. \textit{J. Comb. Theory} 16, 47--56(1974).
		
		
		
		\bibitem{habib-2005} Habib, M., Paul, C.: A simple linear time algorithm for cograph recognition. \textit{Discrete Appl. Mathematics} 145, 183--197(2005).
		
		\bibitem{jung} Jung, H.A.: On a class of posets and the corresponding comparability graphs. \textit{Journal of Combin. Theory} B 24, 125-133(1978)
		
		
		
		
		\bibitem{mdp-14} Maxov\'{a},~J.,  Dubcov\'{a},~M., Pavl\'{i}kov\'{a},~P.,  Turz\'{i}k,~D.: Which k-trees are cover-incomparability graphs. \textit{Discrete Appl.Math.} 167, 222--227(2014).
		
		\bibitem{mpt-09} Maxov\'{a},~J., Pavl\'{i}kov\'{a},~P., Turz\'{i}k,~D.: On the complexity of cover-incomparability graphs of posets. \textit{Order} 26(3), 229--236(2009).

\bibitem{matu-13} Maxov\'{a},~J.,  Turz\'{i}k,~ D.: Which distance-hereditary graphs are cover-incomparability graphs?. \textit{Discrete Appl. Math.} 161, 2095--2100(2013).
		
  
		\bibitem{McConnell-1994} McConnell,~R.M.,  Spinrad,~J.:  Linear-time modular decomposition and efficient transitive orientation of comparability graphs. In \textit{Proceedings of the Fifth Annual ACM-SIAM Symposium on Discrete Algorithms Arlington}, VA, ACM, NewYork, 536--545(1994).
		
		
		
		
		\bibitem{rose_chordal-2}  Rose,~D.J., Tarjan,~R.E., Lueker,~G.S.: Algorithmic Aspects of Vertex Elimination on Graphs. \textit{SIAM Journal on Computing} 5(2), 266-283 (1976).
		
		\bibitem{seinsche} Seinsche, S.: On the property of the class of n-colorable graphs. \textit{J. Comb. Theory (B)}, 191--193(1974).
		
		\bibitem{shibata}  Shibata, Y.: On the tree representation of chordal graphs. \textit{Journal of Graph Theory} 12, 421--428(1988).
		
		
		\bibitem{trajan} Tarjan, R.E.,  Yannakakis, M.:   Simple linear-time algorithms to test chordality of graphs, test acyclicity of hypergraphs, and selectively reduce acyclic hypergraphs.  \textit{SIAM J. Comput.} 13,  566--579(1984).
		
		
		
		
		
	\end{thebibliography}
\end{document}